# Size-biased branching population measures and the multi-type $x \log x$ condition

PETER OLOFSSON

*Mathematics Department, Trinity University, 1 Trinity Place, San Antonio, TX 78212, USA.
E-mail: polofsso@trinity.edu*

We investigate the $x \log x$ condition for a general (Crump–Mode–Jagers) multi-type branching process with a general type space by constructing a size-biased population measure that relates to the ordinary population measure via an intrinsic martingale $W_t$. Sufficiency of the $x \log x$ condition for a non-degenerate limit of $W_t$ is proved and conditions for necessity are investigated.

*Keywords:* general branching process; immigration; size-biased measure; $x \log x$ condition

## 1. Introduction

The $x \log x$ condition is a fundamental concept for supercritical Galton–Watson branching processes, being the necessary and sufficient condition for the process to grow as its mean. In a Galton–Watson process with offspring mean $m = E[X] > 1$, let $Z_n$ denote the number of individuals in the $n$th generation and let $W_n = Z_n/m^n$. Then, $W_n$ is a non-negative martingale and hence $W_n \to W$ for some random variable $W$. The *Kesten–Stigum theorem* is given as follows.

**Theorem 1.1.** *If $E[X \log^+ X] < \infty$, then $E[W] = 1$; if $E[X \log^+ X] = \infty$, then $W = 0$ a.s.*

Here, $\log^+ x = \max(0, \log x)$. It can further be shown that $P(W = 0)$ must either be 0 or equal the extinction probability, hence $E[X \log^+ X] < \infty$ implies that $W > 0$ exactly on the set of non-extinction (see, for example, Athreya and Ney (1972)).

The analog for general single-type branching processes appears in Jagers and Nerman (1984) and a partial result (establishing sufficiency) for general multi-type branching processes appears in Jagers (1989). Lyons, Pemantle and Peres (1995) give a slick proof of the Kesten–Stigum theorem based on comparisons between the Galton–Watson measure and another measure, the *size-biased* Galton–Watson measure, on the space of progeny trees. In Olofsson (1998), these ideas were further developed to analyze general single-type branching processes and the current paper considers general multi-type branching processes with a general type space. In addition to providing a new proof of a known result, size-biased processes also provide tools to further analyze necessity of the $x \log x$ condition.

A crucial concept for the Lyons–Pemantle–Peres (LPP) proof is that of size bias. If the offspring distribution is $\{p_0, p_1, \ldots\}$ and has mean $m = E[X]$, then the size-biased offspring distribution is defined as

$$\widetilde{p}_k = \frac{k p_k}{m}$$

for $k = 0, 1, 2, \ldots$, where we note, in particular, that $\widetilde{p}_0 = 0$. A size-biased Galton–Watson tree is constructed in the following way. Let $\widetilde{X}$ be a random variable that has the size-biased offspring distribution and let the ancestor $v_0$ have a number $\widetilde{X}_0$ of children. Pick one of these at random, call her $v_1$, give her a number $\widetilde{X}_1$ of children and give her siblings ordinary Galton–Watson descendant trees. Pick one of $v_1$'s children at random, call her $v_2$, give her a number $\widetilde{X}_2$ of children, give her sisters ordinary Galton–Watson descendant trees and so on and so forth. With $P_n$ denoting the ordinary Galton–Watson measure restricted to the $n$ first generations, $\widetilde{P}_n$ denoting the measure that arises from the above construction and $W_n = Z_n/m^n$, it can be shown that the relation

$$\mathrm{d}\widetilde{P}_n = W_n \, \mathrm{d}P_n \tag{1.1}$$

holds. Hence, it is the martingale $W_n$ that size-biases the Galton–Watson process. The construction of $\widetilde{P}$ can also be viewed as describing a Galton–Watson process with immigration, where the immigrants are the siblings of the individuals on the path $(v_0, v_1, \ldots)$. Thus, the measure $P$ is the ordinary Galton–Watson measure and the size-biased measure $\widetilde{P}$ is the measure of a Galton–Watson process with immigration, where the i.i.d. immigration group sizes are distributed as $\widetilde{X} - 1$. The relation between $P$ and $\widetilde{P}$ on the space of family trees can now be explored using results for processes with immigration and this provides the final key to the proof.

The general idea of using size-bias in branching processes appeared before LPP. One early example is Joffe and Waugh (1982), where size-biased Galton–Watson processes show up in the study of ancestral trees of randomly sampled individuals. This approach was further explored by Olofsson and Shaw (2002) with a view toward biological applications. An approach similar to LPP appeared in Waymire and Williams (1996), developed simultaneously with, and independently of, LPP. Later applications and extensions of the powerful LPP method include Kurtz *et al.* (1997), Geiger (1999), Athreya (2000), Biggins and Kyprianou (2004) and Lambert (2007).

To make this paper self-contained, we give a short review of general multi-type branching processes and their $x \log x$ condition in the next section. As in the Galton–Watson case, branching processes with immigration are crucial in the proof; for that purpose, we briefly discuss processes with immigration in Section 3, following Olofsson (1996). The size-biased measure on the space of population trees and its relation to the ordinary branching measure is investigated in Section 4 and, in Section 5, sufficiency of the $x \log x$ condition is proved. Finally, in Section 6, we discuss various conditions for necessity.

## 2. The $x \log x$ condition for general branching processes

In a general branching process, individuals are identified by descent. The ancestor is denoted by 0, the children of the ancestor by $1, 2, \ldots$ and so on, so that the individual $x = (x_1, \ldots, x_n)$ is the $x_n$th child of the $x_{n-1}$th child of ... of the $x_1$th child of the ancestor. The set of all individuals can thus be described as

$$I = \bigcup_{n=0}^{\infty} N^n.$$

At birth, each individual is assigned a *type s*, chosen from the *type space S*, equipped with some appropriate $\sigma$-algebra $\mathcal{S}$. The type space can be quite general; usually, it is required to be a complete, separable metric (that is, Polish) space. The type $s$ determines a probability measure $P_s(\cdot)$, the *life law*, on the *life space* $\Omega$, equipped with some appropriate $\sigma$-algebra $\mathcal{F}$. The information provided by a life $\omega \in \Omega$ may differ from one application to another, but it must at least give the *reproduction process* $\xi$ on $S \times R_+$. This process gives the sequence of birth times and types of the children of an individual. More precisely, let $(\tau(k), \sigma(k))$ be random variables on $\Omega$ denoting the birth time (age of the mother) and type of the $k$th child, respectively, and define

$$\xi(A \times [0, t]) = \#\{k : \sigma(k) \in A, \tau(k) \leq t\}$$

for $A \in \mathcal{S}$ and $t \geq 0$. We let $\tau(k) \equiv \infty$ if fewer than $k$ children are born. The *population space* is defined as $\Omega^I$, an outcome $\omega^I$ of which gives the lives of all individuals, together with the $\sigma$-algebra $\mathcal{F}^I$. The set of probability kernels $\{P_s(\cdot), s \in S\}$ defines a probability measure on $(\Omega^I, \mathcal{F}^I)$, the *population measure* $P_s$, where the ancestor's type is $s$.

With each individual $x \in I$, we associate its type $\sigma_x$, its birth time $\tau_x$ and its life $\omega_x$, where $\sigma_x$ is inherited from the mother (a function of the mother's life) and $\omega_x$ is chosen according to the probability distribution $P_{\sigma_x}(\cdot)$ on $(\Omega, \mathcal{F})$. The birth time $\tau_x$ is defined recursively by letting the ancestor be born at time $\tau_0 = 0$ and, if $x$ is the $k$th child of its mother $y$, we let $\tau_x = \tau_y + \tau(k)$. Note that $\tau_x$ and $\tau_y$ denote absolute time, whereas $\tau(k)$ is the mother's age at $x$'s birth.

An important entity is the *reproduction kernel*, defined by

$$\mu(s, \mathrm{d}r \times \mathrm{d}t) = E_s[\xi(\mathrm{d}r \times \mathrm{d}t)],$$

the expectation of $\xi(\mathrm{d}r \times \mathrm{d}t)$ when the mother is of type $s$. This kernel plays the role of $m = E[X]$ in the simple Galton–Watson process and determines the growth rate of the process as $\mathrm{e}^{\alpha t}$, where $\alpha$ is called the *Malthusian parameter*. We assume throughout that the process is supercritical, meaning that $\alpha > 0$. The existence of such an $\alpha$ is not automatic; sufficient conditions may be found in Jagers (1989, 1992). For the rest of this section, we leave out further technical details and assumptions, instead focusing on the main definitions and results. The details can be found in Jagers (1989, 1992) and we simply refer to a process that satisfies all of the conditions needed as a *Malthusian* process.

Given $\mu$ and $\alpha$, we define the kernel $\widehat{\mu}$ as

$$\widehat{\mu}(s, \mathrm{d}r) = \int_0^\infty \mathrm{e}^{-\alpha t} \mu(s, \mathrm{d}r \times \mathrm{d}t) \tag{2.1}$$

which, under certain conditions, has eigenmeasure $\pi$ and eigenfunction $h$ given by

$$\begin{aligned} \pi(\mathrm{d}r) &= \int_S \widehat{\mu}(s, \mathrm{d}r) \pi(\mathrm{d}s), \\ h(s) &= \int_S h(r) \widehat{\mu}(s, \mathrm{d}r), \end{aligned} \tag{2.2}$$

where both $\pi$ and $h\,\mathrm{d}\pi$ can be normed to probability measures. The measure $\pi$ is called the *stable type distribution* and $h(s)$ is called the *reproductive value* of an individual of type $s$. The interpretation of $\pi$ and $h$ is that $\pi$ is the limiting distribution of the type of an individual chosen at random from a population and $h(s)$ is a measure of how reproductive the type $s$ tends to be, in a certain average sense. Moreover, after suitable norming, it can be shown that $h\,\mathrm{d}\pi$ is the probability measure that is the limiting type distribution backward in the family tree from the randomly sampled individual mentioned above. The mean asymptotic age of a random child-bearing in this backward sense is denoted by $\beta$ and satisfies

$$\beta = \int_{S \times S \times R_+} t\mathrm{e}^{-\alpha t} h(r) \mu(s, \mathrm{d}r \times \mathrm{d}t) \pi(\mathrm{d}s) < \infty. \tag{2.3}$$

To count, or measure, the population, *random characteristics* are used. A random characteristic is a real-valued process $\chi$, where $\chi(a)$ gives the contribution to the population of an individual of age $a$. Thus, $\chi$ is a process defined on the life space and by letting $\chi_x$ be the characteristic pertaining to the individual $x$, the $\chi$-counted population is defined as

$$Z_t^\chi = \sum_{x \in I} \chi_x(t - \tau_x),$$

which is the sum of the contributions of all individuals at time $t$ (when the individual $x$ is of age $t - \tau_x$). The simplest example of a random characteristic is $\chi(a) = I_{R_+}(a)$, the indicator for being born, in which case $Z_t^\chi$ is simply the total number of individuals born up to time $t$.

To capture the asymptotics of $Z_t^\chi$, the crucial entity is the *intrinsic martingale* $W_t$, introduced by Nerman (1981) for single-type processes and generalized to multi-type processes in Jagers (1989). For its definition, denote $x$'s mother by $mx$ and let

$$\mathcal{I}_t = \{x : \tau_{mx} \leq t < \tau_x\}, \tag{2.4}$$

the set of individuals whose mothers are born at, or before, time $t$, but who themselves are not yet born at time $t$. The set $\mathcal{I}_t$, sometimes referred to as the "coming generation", generalizes the concept of generation in the Galton–Watson process. Now, let

$$W_t = \frac{1}{h(\sigma_0)} \sum_{x \in \mathcal{I}_t} e^{-\alpha \tau_x} h(\sigma_x), \tag{2.5}$$

the individuals in $\mathcal{I}_t$ summed with time- and type-dependent weights, normed by the reproductive value of the ancestor. It can be shown that $W_t$ is a martingale with respect to the $\sigma$-algebra $\mathcal{F}_t$ generated by the lives of all individuals born before $t$ and that $E_s[W_t] = 1$ for all $s \in S$. Hence, $W_t$ plays the role that $W_n = Z_n/m^n$ does in the Galton–Watson process and the limit of $Z_t^\chi$ turns out to involve the martingale limit $W = \lim_{t \to \infty} W_t$. The main convergence result is of the form

$$e^{-\alpha t} Z_t^\chi \to \frac{E_\pi[\widehat{\chi}(\alpha)]}{\alpha \beta} h(s) W$$

$P_s$-almost surely for $\pi$-almost all $s \in S$ as $t \to \infty$. Here, $\sigma_0 = s$ is the type of the ancestor, $E_\pi[\cdot] = \int_S E_s[\cdot] \pi(\mathrm{d}s)$ and $\widehat{\chi}(\alpha)$ is the Laplace transform of $\chi(a)$ evaluated at the point $\alpha$. As in the Galton–Watson case, the question is when the martingale limit $W$ is non-degenerate. As $W_t \to W$ $P_s$-a.s. and $E_s[W_t] = 1$, $L^1$-convergence with respect to $P_s$ is equivalent to $E_s[W] = 1$ (Durrett (2005), page 258). Note that although it is the process $Z_t^\chi$ that is of interest and not $W_t$ itself, the asymptotics are determined by $W_t$, one of many examples in probability of the usefulness of finding an embedded martingale.

We are ready to formulate the general $x\log x$ condition and the main convergence result. For the reproduction process $\xi$, define the transform

$$\bar{\xi} = \int_{S \times R_+} e^{-\alpha t} h(r) \xi(\mathrm{d}r \times \mathrm{d}t) \tag{2.6}$$

which plays the role of $X$ in the Galton–Watson process (in fact, in that case, $\bar{\xi} = X/m$). For future reference, let us also state an alternative representation of $\bar{\xi}$. Denote the sequence of birth times and types in the process $\xi$ by $\tau(1), \sigma(1), \tau(2), \sigma(2), \ldots$ and so on. Then,

$$\bar{\xi} = \sum_{i=1}^{\infty} e^{-\alpha \tau(i)} h(\sigma(i)). \tag{2.7}$$

The $x\log x$ condition and convergence result are given in the following theorem from Jagers (1989).

**Theorem 2.1.** *Consider a general multi-type supercritical Malthusian branching process with*

$$E_\pi[\bar{\xi} \log^+ \bar{\xi}] < \infty.$$

*Then, $E_s[W] = 1$ for $\pi$-almost all $s$, from which it follows that*

$$e^{-\alpha t} Z_t^\chi \to \frac{E_\pi[\widehat{\chi}(\alpha)]}{\alpha \beta} h(s) W$$

in $L^1(P_s)$ for $\pi$-almost all $s$.

## 3. Processes with immigration

As mentioned in the Introduction, branching processes with immigration are crucial to our proof and in this section, we state the main result for such processes. Consider a general branching process where new individuals immigrate into the population according to some point process $\eta(\mathrm{d}r \times \mathrm{d}t)$ with points of occurrence and types $(\tau_1, \sigma_1), (\tau_2, \sigma_2), \ldots$. The $k$th immigrant initiates a branching process according to the population measure $P_{\sigma_k}$. The immigration process has the transform

$$\bar\eta = \int_0^\infty \mathrm{e}^{-\alpha t} h(r) \eta(\mathrm{d}r \times \mathrm{d}t) = \sum_{k=1}^\infty \mathrm{e}^{-\alpha \tau_k} h(\sigma_k)$$

and it can be shown that the process $W_t$ is now a submartingale rather than a martingale (which is intuitively clear because offspring of immigrants may be added to the set $\mathcal{I}_t$). The limit of $W_t$ is therefore not automatically finite, but needs a condition on the immigration process, established by the following lemma from Olofsson (1996).

**Lemma 3.1.** *If $\bar\eta < \infty$ a.s., then $W_t \to W$ a.s. as $t \to \infty$, where $W < \infty$ a.s.*

## 4. The size-biased population measure

Recall that the LPP size-biased Galton–Watson measure was constructed from the size-biased offspring distribution. General branching processes require a more general concept of size-bias. In a general process, the offspring random variable $X$ is replaced by the reproduction process $\xi$, the size of which is properly measured by the transform $\bar\xi$ which leads to the following definition.

**Definition 4.1.** *The size-biased life law $\widetilde P_s$ is defined as*

$$\widetilde P_s(\mathrm{d}\omega) = \frac{\bar\xi(\omega)}{h(s)} P_s(\mathrm{d}\omega).$$

From Jagers (1992), we know that the eigenfunction $h$ is finite and strictly positive, so $\widetilde P$ is well defined. The following lemma follows immediately from the definition of $\widetilde P_s$.

**Lemma 4.2.** *Let $P_s$ and $\widetilde P_s$ be as above and denote the set of realizations of reproduction processes by $\Gamma$, equipped with a $\sigma$-algebra $\mathcal{G}$. Then,*

(i) *for $A \in \mathcal{F}$,*

$$\widetilde P_s(A) = \frac{E_s[\bar\xi; A]}{h(s)};$$

(ii) *for every $\mathcal{G}$-measurable function $g\colon \Gamma \to R$,*

$$\widetilde{E}_s[g(\xi)] = \frac{E_s[\bar{\xi}g(\xi)]}{h(s)}.$$

Note that $\widetilde{P}_s$ is indeed a probability measure for all $s \in S$ because

$$\widetilde{P}_s(\Omega) = \frac{1}{h(s)} E_s[\bar{\xi}] = 1,$$

where $E_s[\bar{\xi}] = h(s)$ follows from the definition of $\bar{\xi}$ in (2.6), together with (2.1) and (2.2). Also, note that a size-biased reproduction process always contains points because

$$\widetilde{P}_s(\xi(S \times R_+) = 0) = \frac{1}{h(s)} E_s[\bar{\xi}; \xi(S \times R_+) = 0] = 0,$$

in analogy with the size-biased offspring distribution in a Galton–Watson process ($\xi(S \times R_+)$ is the total number of offspring of an individual).

To construct the size-biased population measure, let $\widetilde{P}_s^t$ and $P_s^t$ denote the restrictions of the measures $P_s$ and $\widetilde{P}_s$ to the $\sigma$-algebra $\mathcal{F}_t$. The goal is to construct a measure $\widetilde{P}_s$ on $(\Omega^I, \mathcal{F}^I)$ that is such that

$$\widetilde{P}_s^t(\mathrm{d}\omega^I) = W_t(\omega^I) P_s^t(\mathrm{d}\omega^I)$$

for all $t$, where $W_t$ is the intrinsic martingale defined in (2.5). This measure is the direct extension of the size-biased measures from Lyons *et al.* (1995) and Olofsson (1998). The construction also involves the set $\mathcal{I}_t$, defined in (2.4), whose individuals all have mothers that are born up to time $t$. Thus, the type and birth time of an individual in $\mathcal{I}_t$ is measurable with respect to $\mathcal{F}_t$, which implies that $W_t$ is also measurable with respect to $\mathcal{F}_t$.

The construction of the size-biased population measure extends the construction in Olofsson (1998) as follows. Start with the ancestor, now called $v_0$, and choose her life $\omega_0$ according to the size-biased distribution $\widetilde{P}_s(\mathrm{d}\omega_0) = \frac{\bar{\xi}_0}{h(s)} P_s(\mathrm{d}\omega_0)$. Pick one of her children, born in the reproduction process $\xi_0$, such that the $i$th child is chosen with probability $\frac{e^{-\alpha \tau_i} h(\sigma_i)}{\bar{\xi}_0}$. Call this child $v_1$, let her start a population according to the size-biased population law $\widetilde{P}_{\sigma_i}$ and give her sisters independent descendant trees such that sister $j$ initiates a branching process according to the regular population law $P_{\sigma_j}$. Continue in this way and define the measure $\widetilde{P}_s$ to be the joint distribution of the random tree and the random path $(v_0, v_1, \ldots)$. We shall borrow a term from Athreya (2000) and refer to the path $(v_0, v_1, \ldots)$ of chosen individuals as the *spine*.

Now, fix an individual $x$ in the set $\mathcal{I}_t$ defined in (2.4) and consider the probability $\widetilde{P}_s$, constrained by the individual $x$ being chosen to be in the spine. Specifically, if $S(x,t)$ denotes the event that the individual $x$ in $\mathcal{F}_t$ is chosen to be in the spine and $A \in \mathcal{F}^I$, then we consider the measure $\widetilde{P}_s(\cdot\,;x)$ defined by

$$\widetilde{P}_s^t(A;x) = \widetilde{P}_s^t(A \cap S(x,t)). \tag{4.1}$$

Denote by $i$ the individual in the first generation from whom $x$ stems, that is, $x = (i,y)$ for some $y$. Hence, if $x$ is in the $n$th generation, then it is of the form $x = (x_1, x_2, \ldots, x_n)$, where $x_1 = i$, and we let $y = (x_2, \ldots, x_n)$. In words, $y$ is the same individual as $x$ when $i$ is viewed as the ancestor. Let $\omega^{(j)}$ denote the lives of all individuals when $j$ is viewed as the ancestor to obtain

$$\widetilde{P}_s^t(\mathrm{d}\omega^I;x) = \frac{\bar{\xi}_0}{h(s)} P_s(\mathrm{d}\omega_0) \cdot \frac{\mathrm{e}^{-\alpha\tau_i} h(\sigma_i)}{\bar{\xi}_0} \cdot \widetilde{P}_{\sigma_i}^{t-\tau_i}(\mathrm{d}\omega^{(i)};y) \cdot \prod_{j \neq i} P_{\sigma_j}^{t-\tau_j}(\mathrm{d}\omega^{(j)}), \tag{4.2}$$

where the first factor describes the size-biased choice of life of the ancestor. The second factor is the probability that the individual $i$ in the first generation is chosen to be in the spine and the third factor describes the size-biased probability measure of the process starting from $i$, constrained by the individual $y$ being in the spine. Finally, the fourth factor describes the regular population measures stemming from the individuals in the first generation who are not chosen to be in the spine.

The following proposition states the desired relation between the size-biased measure $\widetilde{P}_s^t$ and the regular population measure $P_s^t$.

**Proposition 4.3.** *Let $\widetilde{P}_s^t$ and $P_s^t$ be the restrictions of $\widetilde{P}_s$ and $P_s$ to the $\sigma$-algebra $\mathcal{F}_t$ and let $W_t$ be as in (2.5). Then,*

$$\frac{\mathrm{d}\widetilde{P}_s^t}{\mathrm{d}P_s^t} = W_t.$$

**Proof.** Let $\widetilde{P}_s(\cdot\,;x)$ be as in (4.1). Then,

$$\widetilde{P}_s^t(\mathrm{d}\omega^I;x) = \frac{\mathrm{e}^{-\alpha\tau_x} h(\sigma_x)}{h(s)} P_s^t(\mathrm{d}\omega^I).$$

Further, note that for the regular population measure, we have

$$P_s^t(\mathrm{d}\omega^I) = P_s(\mathrm{d}\omega_0) \prod_{j=1}^{\xi_0(t)} P_{\sigma_j}^{t-\tau_j}(\mathrm{d}\omega^{(j)})$$

$$= P_s(\mathrm{d}\omega_0) P_{\sigma_i}^{t-\tau_i}(\mathrm{d}\omega^{(i)}) \prod_{j \neq i} P_{\sigma_j}^{t-\tau_j}(\mathrm{d}\omega^{(j)})$$

$$= \bar{\xi}_0 P_s(\mathrm{d}\omega_0) \frac{h(\sigma_i)}{\bar{\xi}_0} \frac{1}{h(\sigma_i)} P_{\sigma_i}^{t-\tau_i}(\mathrm{d}\omega^{(i)}) \prod_{j \neq i} P_{\sigma_j}^{t-\tau_j}(\mathrm{d}\omega^{(j)}).$$

Now, let $\tau_y(i)$ and $\sigma_y(i)$ denote the birth time and type of the individual $y$ when $i$ is viewed as the ancestor. We then have $\tau_x = \tau_i + \tau_y(i)$ and $\sigma_x = \sigma_y(i)$. Multiply $P_s^t(\mathrm{d}\omega^I)$ by $\frac{\mathrm{e}^{-\alpha\tau_x}h(\sigma_x)}{h(s)}$ to obtain

$$\frac{\mathrm{e}^{-\alpha\tau_x}h(\sigma_x)}{h(s)}P_s^t(\mathrm{d}\omega^I)$$

$$= \frac{\bar{\xi}_0}{h(s)}P_s(\mathrm{d}\omega_0) \cdot \frac{\mathrm{e}^{-\alpha\tau_i}h(\sigma_i)}{\bar{\xi}_0} \cdot \frac{\mathrm{e}^{-\alpha\tau_y(i)}h(\sigma_y(i))}{h(\sigma_i)} \cdot P_{\sigma_i}^{t-\tau_i}(\mathrm{d}\omega^{(i)}) \cdot \prod_{j \neq i} P_{\sigma_j}^{t-\tau_j}(\mathrm{d}\omega^{(j)})$$

$$= \frac{\bar{\xi}_0}{h(s)}P_s(\mathrm{d}\omega_0) \cdot \frac{\mathrm{e}^{-\alpha\tau_i}h(\sigma_i)}{\bar{\xi}_0} \cdot \widetilde{P}_{\sigma_i}^{t-\tau_i}(\mathrm{d}\omega^{(i)};y) \cdot \prod_{j \neq i} P_{\sigma_j}^{t-\tau_j}(\mathrm{d}\omega^{(j)})$$

$$= \widetilde{P}_s^t(\mathrm{d}\omega^I;x),$$

by (4.2). Finally, sum over $x \in \mathcal{I}_t$ to get

$$\widetilde{P}_s^t(\mathrm{d}\omega^I) = \sum_{x \in \mathcal{I}_t} \frac{\mathrm{e}^{-\alpha\tau_x}h(\sigma_x)}{h(s)}P_s^t(\mathrm{d}\omega^I) = W_t P_s^t(\mathrm{d}\omega^I),$$

where $\omega^I$ is suppressed, but understood as the argument of $\tau_x$, $\sigma_x$ and $W_t$.    □

To summarize, restricted to the $\sigma$-algebra $\mathcal{F}_t$, the size-biased population measure $\widetilde{P}_s$ relates to the regular population measure $P_s$ via the Radon–Nikodym derivative $W_t = \frac{\mathrm{d}\widetilde{P}_s^t}{\mathrm{d}P_s^t}$, a relation that is a straightforward extension from Olofsson (1998), which, in turn, is the straightforward extension of the original LPP method. This result is also similar in nature to Proposition 1 in Athreya (2000), which deals with a different martingale.

The individuals $v_0, v_1, \ldots$ in the spine are of particular interest in our analysis. From now on, we will use $\sigma_0, \sigma_1, \ldots$ to denote the types of the individuals in the spine. The inter-arrival times are denoted $T_1, T_2, \ldots$ that is, $T_k$ is the time between the appearances of the $(k-1)$th and $k$th individual. In the single-type case treated in Olofsson (1998), the individuals in the spine have lives that are i.i.d., but the situation is now much different, with dependence between consecutive individuals introduced via types. As we shall see later, this dependence is also what prevents the proof of necessity of the $x \log x$ condition from carrying over from the single-type case.

We now state important properties of the sequences of types and inter-arrival times in two lemmas. The first lemma deals solely with the type sequence. For the rest of this section, we use the notation $\widetilde{P}$ rather than $\widetilde{P}_s$ since the conditional probabilities we consider do not depend on the initial type.

**Lemma 4.4.** *The sequence of types $(\sigma_0, \sigma_1, \ldots)$ in the spine is a homogeneous Markov chain with transition probabilities*

$$\widetilde{P}(\sigma_{k+1} \in \mathrm{d}r | \sigma_k = s) = \frac{h(r)}{h(s)}\widehat{\mu}(s, \mathrm{d}r)$$

*and stationary distribution $\nu(\mathrm{d}s) = h(s)\pi(\mathrm{d}s)$.*

**Proof.** In a generic reproduction process $\xi$, denote the birth time and type of the $i$th offspring by $\tau(i)$ and $\sigma(i)$, respectively. Note the difference between $\sigma(i)$ and $\sigma_i$, the latter being the type of the $i$th individual in the spine. The transition probabilities satisfy

$$\widetilde{P}(\sigma_{k+1} \in \mathrm{d}r | \sigma_k = s)$$
$$= \sum_i \widetilde{P}(\sigma_{k+1} \in \mathrm{d}r, v_{k+1} = i | \sigma_k = s) = \sum_i \widetilde{E}_s\left[\frac{\mathrm{e}^{-\alpha\tau(i)}h(\sigma(i))}{\bar{\xi}}\delta_{\sigma(i)}(\mathrm{d}r)\right]$$
$$= \frac{1}{h(s)}\sum_i E_s[\mathrm{e}^{-\alpha\tau(i)}h(\sigma(i))\delta_{\sigma(i)}(\mathrm{d}r)] = \frac{1}{h(s)}\int_0^\infty E_s[\mathrm{e}^{-\alpha t}h(r)\xi(\mathrm{d}r \times \mathrm{d}t)]$$
$$= \frac{h(r)}{h(s)}\widehat{\mu}(s, \mathrm{d}r),$$

where we have used Lemma 4.2 applied to the function

$$g(\xi) = \sum_i \frac{\mathrm{e}^{-\alpha\tau_i}h(\sigma_i)}{\bar{\xi}}.$$

Next, let $\nu(\mathrm{d}s) = h(s)\pi(\mathrm{d}s)$. As

$$\int_S \widehat{\mu}(s, \mathrm{d}r)\pi(\mathrm{d}s) = \pi(\mathrm{d}r),$$

we get

$$\int_{s \in S} \frac{h(r)}{h(s)}\widehat{\mu}(s, \mathrm{d}r)\nu(\mathrm{d}s) = \nu(\mathrm{d}r)$$

and thus the Markov chain of types in the spine has stationary distribution $\nu = h \, \mathrm{d}\pi$. $\square$

The second lemma deals with the sequence of types and inter-arrival times of the individuals in the spine.

**Lemma 4.5.** *The sequence of types and inter-arrival times $(\sigma_0, T_1, \sigma_1, T_2, \ldots)$ of the individuals in the spine constitutes a Markov renewal process with transition kernel*

$$\widetilde{P}(T_{k+1} \in \mathrm{d}t, \sigma_{k+1} \in \mathrm{d}r | \sigma_k = s) = \frac{h(r)}{h(s)}\mathrm{e}^{-\alpha t}\mu(s, \mathrm{d}r \times \mathrm{d}t)$$

*and the expected value of $T_k$ when $\sigma_0 \sim \nu$ is*

$$\widetilde{E}_\nu[T_k] = \beta < \infty,$$

*where $\beta$ was defined in (2.3).*

**Proof.** Similarly to the proof of Lemma 4.4, we get

$$\widetilde{P}(T_{k+1} \in \mathrm{d}t, \sigma_{k+1} \in \mathrm{d}r | \sigma_k = s)$$
$$= \sum_i \widetilde{P}(T_{k+1} \in \mathrm{d}t, \sigma_{k+1} \in \mathrm{d}r, v_{k+1} = i | \sigma_k = s)$$
$$= \frac{1}{h(s)} E_s[\mathrm{e}^{-\alpha t} h(r) \xi(\mathrm{d}r \times \mathrm{d}t)]$$
$$= \frac{h(r)}{h(s)} \mathrm{e}^{-\alpha t} \mu(s, \mathrm{d}r \times \mathrm{d}t)$$

and the expected value of $T_k$ when $\sigma_0$ is chosen according to the stationary distribution $\nu = h \, \mathrm{d}\pi$ is

$$\widetilde{E}_\nu[T_k] = \int_{S \times [0, \infty)} t \frac{h(r)}{h(s)} \mathrm{e}^{-\alpha t} \mu(s, \mathrm{d}r \times \mathrm{d}t) h(s) \pi(\mathrm{d}s)$$
$$= \int_{S \times [0, \infty)} t \mathrm{e}^{-\alpha t} h(r) \mu(s, \mathrm{d}r \times \mathrm{d}t) \pi(\mathrm{d}s) = \beta,$$

by (2.3). □

There is an interesting connection between the size-biased measure and the *stable population measure* from Jagers (1992). The latter is an asymptotic probability measure that is centered around a randomly sampled individual as $t \to \infty$. In such a stable population, the randomly sampled individual is born in a point process that has the size-biased distribution, the asymptotic type distribution as time goes backward through the individual's line of descent is $h \, \mathrm{d}\pi$ and the asymptotic mean age at a random child-bearing is $\beta$. The transition probabilities in this backward chain also involve $\widehat{\mu}(s, \mathrm{d}r)$, but have weights that are expressed in terms of $\pi$ rather than $h$, as we have in the size-biased measure where time goes forward. This relation becomes clearer in a finite-type Galton–Watson process where $\pi$ and $h$ are simply the left and right eigenvectors of the mean reproduction matrix.

## 5. Sufficiency of the $x \log x$ condition

We will soon be ready to prove the general $x \log x$ theorem, Theorem 2.1, the key to which is the relation between $\widetilde{P}_s$ and $P_s$. Recall that the two are related through $W_t$, which is a martingale under $P_s$ and a submartingale under $\widetilde{P}_s$. The following lemma relates the limiting behavior of $W_t$ under $P_\pi(\cdot) = \int_S P_s(\cdot) \pi(\mathrm{d}s)$ to its limiting behavior under $\widetilde{P}_\nu(\cdot) = \int_S \widetilde{P}_s(\cdot) h(s) \pi(\mathrm{d}s)$.

**Lemma 5.1.** *Let $W = \limsup_t W_t$. Then,*

(i) $\widetilde{P}_\nu(W = \infty) = 0 \Rightarrow E_\pi[W] = 1$;
(ii) $\widetilde{P}_\nu(W = \infty) = 1 \Rightarrow E_\pi[W] = 0$.

**Proof.** By Durrett (2005), page 239,

$$\widetilde{P}_\nu(A) = \widetilde{E}_\nu[W; A] + \widetilde{P}_\nu(A \cap \{W = \infty\}).$$

If $\widetilde{P}_\nu(W = \infty) = 0$, then we have $\widetilde{P}_\nu(A) = \widetilde{E}_\nu[W; A]$ and get

$$\widetilde{P}_\nu(W = 0) = \widetilde{E}_\nu[W; W = 0] = 0.$$

Hence,

$$\widetilde{E}_\nu[W] = \widetilde{E}_\nu[W : W > 0] = \widetilde{P}_\nu(W > 0) = 1.$$

Moreover, as

$$1 = \widetilde{E}_\nu[W] = \int_S \widetilde{E}_s[W] \nu(\mathrm{d}s)$$

and as Fatou's lemma implies that $\widetilde{E}_s[W] \leq 1$ for all $s$, we must have $\widetilde{E}_s[W] = 1$ for $\nu$ almost all $s \in S$. As $\pi \ll \nu$, we also get $E_\pi[W] = 1$.

Next, suppose that $\widetilde{P}_\nu(W = \infty) = 1$. As $W$ is an a.s. finite martingale limit under $P_s$, for $\pi$-almost all $s \in S$, we have

$$P_\nu(W = \infty) = \int_S P_s(W = \infty) h(s) \pi(\mathrm{d}s) = 0.$$

Hence, the measures $\widetilde{P}_\nu(\cdot)$ and $E_\nu[W; \cdot]$ are mutually singular and, as $\widetilde{P}_\nu(W > 0) > 0$, we get

$$E_\nu[W] = E_\nu[W; W > 0] = 0,$$

which implies that $E_\pi[W] = 0$ as well. The proof of the lemma is thus complete. □

Another connection between expected values under $P_\pi$ and $\widetilde{P}_\nu$, directly involving the $x \log x$ condition, is given by the following corollary to Lemma 4.2.

**Corollary 5.2.**

$$E_\pi[\bar{\xi} \log^+ \bar{\xi}] = \widetilde{E}_\nu[\log^+ \bar{\xi}].$$

**Proof.** Choose $g(\xi) = \log^+ \bar{\xi}$ in Lemma 4.2 to obtain

$$E_\pi[\bar{\xi} \log^+ \bar{\xi}] = \int_S E_s[\bar{\xi} \log^+ \bar{\xi}] \pi(\mathrm{d}s)$$

$$= \int_S h(s) \widetilde{E}_s[\log^+ \bar{\xi}] \pi(\mathrm{d}s)$$

$$= \widetilde{E}_\nu[\log^+ \bar{\xi}]. \qquad \square$$

The logical structure of the proof of Theorem 2.1 is as follows. Assume that $E_\pi[\xi \log^+ \xi] < \infty$. By Corollary 5.2, we then have $\widetilde{E}_\nu[\log^+ \bar{\xi}] < \infty$. If we can show that this, in turn, implies that $W < \infty$ almost surely with respect to $\widetilde{P}_\nu$, then we can invoke Lemma 5.1(i) to conclude that $E_\pi[W] = 1$, after which the proof is more or less complete. The gap is filled by the next lemma, which utilizes results for general branching processes with immigration. It has been mentioned that such processes play a vital role in the LPP method and, in the current setting, we observe that the individuals off the spine constitute a general branching process with immigration, the immigrants being the siblings of the individuals $v_1, v_2, \ldots$ in the spine in an obvious extension of the Galton–Watson case. To describe the immigration process formally, let $I_{j,k}$ be the indicator of the event that $v_{j-1}$'s $k$th child is *not* chosen to become $v_j$, denote the immigration time of the $j$th immigrant by $\tau_j$ and denote the birth time and type of the $k$th individual in $\xi_j$ by $\tau_k(j)$ and $\sigma_k(j)$, respectively. The immigration process $\eta$ is

$$\eta(\mathrm{d}s \times \mathrm{d}t) = \sum_{j,k} \delta_{\sigma_k(j)}(\mathrm{d}s) \delta_{\tau_k(j)}(\mathrm{d}t - \tau_j) I_{j,k},$$

which has

$$\bar{\eta} = \sum_{j,k} h(\sigma_k(j)) \mathrm{e}^{-\alpha \tau_j} \mathrm{e}^{-\alpha \tau_k(j)} I_{j,k}. \tag{5.1}$$

We are now ready for the last lemma needed for the proof of Theorem 2.1.

**Lemma 5.3.** *Consider a general branching process with immigration process $\eta$ as above. If $\widetilde{E}_\nu[\log^+ \bar{\xi}] < \infty$, then $W = \lim_{t \to \infty} W_t$ exists and is finite $\widetilde{P}_\nu$-a.s.*

**Proof.** First, note that

$$\widetilde{E}_\nu[\log^+ \bar{\xi}] < \infty \quad \Rightarrow \quad \sum_n \widetilde{P}_\nu(\log^+ \bar{\xi} > cn) < \infty \qquad \text{for all } c > 0.$$

Now, consider the sequence $\bar{\xi}_0, \bar{\xi}_1, \ldots$ for the individuals in the spine. Recalling that $\nu$ is the stationary distribution of the Markov chain of types in the spine, we conclude that $\widetilde{P}_\nu(\log^+ \bar{\xi} > cn) = \widetilde{P}_\nu(\log^+ \bar{\xi}_n > cn)$ for all $n$, which gives

$$\sum_n \widetilde{P}_\nu(\log^+ \bar{\xi}_n > cn) < \infty.$$

By the first Borel–Cantelli lemma, we get

$$\widetilde{P}_\nu(\log^+ \bar{\xi}_n > cn \text{ i.o.}) = 0, \tag{5.2}$$

which, by (5.1), gives

$$\bar{\eta} \leq \sum_{j,k} h(\sigma_k(j)) \mathrm{e}^{-\alpha \tau_j} \mathrm{e}^{-\alpha \tau_k(j)}$$

$$= \sum_{j=1}^{\infty} \mathrm{e}^{-\alpha \tau_j} \bar{\xi}_j < \infty \qquad \widetilde{P}_\nu\text{-a.s.}$$

since the $\tau_j$ are sums of the $T_k$, which, being the regeneration times in a Markov renewal process, obey the strong law of large numbers (Alsmeyer (1994)) so that $\tau_j \sim \beta j$ almost surely as $j \to \infty$ (recall that $\widetilde{E}_\nu[T_k] = \beta < \infty$). The last sum above is a.s. finite because of the subexponential growth of the $\bar{\xi}_n$ established in (5.2). By Lemma 3.1, we conclude that $\lim_t W_t$ exists and is finite $\widetilde{P}_\nu$-a.s. □

We now have in place all of the preliminaries needed to prove Theorem 2.1.

**Proof of Theorem 2.1.** As $E_\pi[\bar{\xi} \log^+ \bar{\xi}] = \widetilde{E}_\nu[\log^+ \bar{\xi}] < \infty$, Lemma 5.3(i) implies that $\widetilde{P}_\nu(W = \infty) = 0$ and Lemma 5.1(ii) gives $E_\pi[W] = 1$. Moreover, as

$$E_\pi[W] = \int_S E_s[W] \pi(\mathrm{d}s)$$

and as Fatou's lemma implies that $E_s[W] \leq 1$ for all $s$, we must have $E_s[W] = 1$ for $\pi$-almost all $s \in S$. The proof is thus complete. □

## 6. Necessity of the $x \log x$ condition

For single-type processes, the condition of having a finite $x \log x$ moment is both sufficient and necessary. Using the size-bias method, this can be established by using the first and second Borel–Cantelli lemmas, respectively. However, in the multi-type setting, the main result in Jagers (1989) establishes only sufficiency and it is currently not known whether necessity holds. The method of size-biased branching processes provides a way of investigating necessity and although the second Borel–Cantelli lemma cannot be used due to dependence, more general versions can be employed. This section is exploratory in nature and does not provide any definite solutions, but rather outlines two different approaches to establish necessity under additional conditions through the conditional Borel–Cantelli lemma and the Kochen–Stone lemma, respectively.

By "the $x \log x$ condition", we mean the condition that $E_\pi[\bar{\xi} \log^+ \bar{\xi}] < \infty$. Hence, to establish necessity, we need to assume that $E_\pi[\bar{\xi} \log^+ \bar{\xi}] = \infty$ and show that this assumption implies that $E_\pi[W] = 0$, invoking Lemma 5.1(ii) in an intermediate step. The logical structure parallels that of the proof of sufficiency: if $E_\pi[\bar{\xi} \log^+ \bar{\xi}] = \infty$, then Lemma 5.2 yields that $\widetilde{E}_\nu[\log^+ \bar{\xi}] = \infty$, which implies that

$$\sum_n \widetilde{P}_\nu(\log^+ \bar{\xi} > cn) = \infty.$$

Since $\nu$ is the stationary distribution for the types in the spine, this further implies that

$$\sum_n \widetilde{P}_\nu(\log^+ \bar{\xi}_n > cn) = \infty.$$

In the proof of sufficiency, this sum was finite and the first Borel–Cantelli lemma could be invoked to conclude that $\widetilde{P}_\nu(\log^+ \bar{\xi}_n > cn \text{ i.o.}) = 0$, leading to the rest of the proof. However, as the reproduction process $\xi_{n+1}$ is chosen according to a probability distribution that depends on the type $\sigma_n$ which is determined by the parent reproduction process $\xi_n$, we cannot assume that the $\bar{\xi}_n$ are independent. Hence, we cannot invoke the second Borel–Cantelli lemma to conclude that, almost surely, $\log^+ \bar{\xi}_n > cn$ i.o. This constitutes a difference from the single-type case treated in Olofsson (1998), where sufficiency and necessity were established using the first and second Borel–Cantelli lemmas, respectively. Below, we establish necessity of the $x \log x$ condition under various additional conditions.

The conditional Borel–Cantelli lemma states that if $\{\mathcal{F}_n\}$ is a filtration and $\{A_n\}$ a sequence of events with $A_n \in \mathcal{F}_n$, then

$$\{A_n \text{ i.o.}\} = \left\{\sum_{n=1}^\infty P(A_n | \mathcal{F}_{n-1}) = \infty\right\};$$

see Durrett (2005). For us, $A_n = \{\log^+ \bar{\xi}_n > cn\}$, the $\sigma$-algebra $\mathcal{F}_{n-1}$ gives the type $\sigma_n$ of the $n$th individual in the spine and we get

$$\widetilde{P}_s(\log^+ \bar{\xi}_n > cn | \mathcal{F}_{n-1}) = \widetilde{P}_{\sigma_n}(\log^+ \bar{\xi} > cn).$$

The question then becomes under which conditions

$$\sum_{n=1}^\infty \widetilde{P}_{\sigma_n}(\log^+ \bar{\xi} > cn) = \infty \qquad \widetilde{P}_\nu\text{-a.s.}$$

given that

$$\sum_{n=1}^\infty \widetilde{P}_\nu(\log^+ \bar{\xi} > cn) = \infty.$$

One more step is necessary in order to invoke part (ii) of Lemma 5.1, namely, to argue that $\log^+ \bar{\xi}_n > cn$ i.o. implies that $W = \limsup_t W_t = \infty$ $\widetilde{P}_\nu$-a.s. We state this as a lemma.

**Lemma 6.1.** *If $\widetilde{P}_\nu(\log^+ \bar{\xi}_n > cn \text{ i.o.}) = 1$, then $\widetilde{P}_\nu(W = \infty) = 1$.*

**Proof.** Consider $W_{\tau_n}$, the value of $W_t$ at the time of the arrival of the $n$th immigrant $v_n$. All of the children of this immigrant except the one chosen to become $v_{n+1}$ belong to $\mathcal{I}_{\tau_n}$, so these children form a subset of $\mathcal{I}_{\tau_n}$. The $k$th child is born at time $\tau_n + \tau(k)$ and has type $\sigma(k)$, where $\tau(k)$ and $\sigma(k)$ are the points in the reproduction process $\xi_n$ of $v_n$. Let $I_k$ be the indicator of the event that the $k$th child of $v_n$ is *not* chosen to become $v_{n+1}$ to obtain

$$W_{\tau_n} = \frac{1}{h(\sigma_0)} \sum_{x \in \mathcal{I}_{\tau_n}} h(\sigma_x) e^{-\alpha \tau_x}$$

$$\geq \frac{e^{-\tau_n}}{h(\sigma_0)} \sum_{k=1}^{\infty} h(\sigma(k)) e^{-\alpha \tau(k)} I_k$$

$$\geq \frac{e^{-\tau_n}}{h(\sigma_0)} (\bar{\xi}_n - C),$$

where $C = \sup_s h(s) < \infty$ (Jagers (1989)) and we recall (2.7). Adjust the last part of the proof of Lemma 5.3 to conclude that $W = \limsup_t W_t = \infty$ $\widetilde{P}_\nu$-a.s. □

Our first result establishes necessity of the $x \log x$ condition under the additional assumption that there is a type that is revisited infinitely often in the spine.

**Proposition 6.2.** *If the Markov chain of types $\sigma_0, \sigma_1, \ldots$ in the spine has one positive recurrent state $r$ such that $E_r[\bar{\xi} \log^+ \bar{\xi}] = \infty$ and if, for $\pi$-almost all starting types $s$, there exists $k$ such that $\sigma_k = r$ a.s., then the $x \log x$ condition is necessary.*

**Note.** *The Markov chain is on a general state space, but here we use "positive recurrence" in its elementary meaning, namely, that $\widetilde{E}_r[T_r] < \infty$, where $T_r = \inf\{n > 0 : \sigma_n = r\}$.*

**Proof of Proposition 6.2.** The result follows from the following observation regarding infinite series. Let $a_n \geq 0$ be a decreasing sequence of real numbers such that $\sum_n a_n = \infty$, let $X_1, X_2, \ldots$ be i.i.d. non-negative and integer-valued random variables with finite mean $\mu$ and let $T_n = X_1 + X_2 + \cdots + X_n$. Then,

$$\sum_{n=1}^{\infty} a_{T_n} = \infty \quad \text{a.s.}$$

This holds because if $k \geq \mu$ is an integer, then $a_{T_n} \geq a_{nk}$ a.s. for large $n$, by the strong law of large numbers, and, obviously, $\sum_n a_{nk} = \infty$ for all fixed $k$. Finally, apply this result to $a_n = \widetilde{P}_s(\log^+ \bar{\xi} > cn)$ with $X_1, X_2, \ldots$ being the consecutive inter-return times to the state $s$ in the Markov chain of types in the spine (that is, the $X_k$ are i.i.d. copies of $T_r$ above). □

One instance where the assumption of Proposition 6.2 follows from the $x \log x$ condition is if the type space is finite. Indeed, if there are $n$ types, we have

$$\infty = E_\pi[\bar{\xi} \log^+ \bar{\xi}] = \sum_{k=1}^{n} E_r[\bar{\xi} \log^+ \bar{\xi}] \pi(r),$$

which implies that we must have $E_r[\bar{\xi} \log^+ \xi] = \infty$ for at least one $r$. In particular, this establishes necessity of the $x \log x$ condition for the ordinary multi-type Galton–Watson process with a finite type space (under the usual assumptions of positive regularity and non-singularity; see Athreya and Ney (1972)).

Another approach is to consider the rate of convergence of the type chain $\sigma_0, \sigma_1, \ldots$ toward its stationary distribution $\nu$. To simplify the analysis, let $Y = [\log^+ \bar{\xi}]$, the integer part of $\log^+ \bar{\xi}$, which has finite mean under $\widetilde{P}_\nu$ if and only if $\log^+ \bar{\xi}$ does. We then have the following.

**Proposition 6.3.** *Suppose that*

$$\sum_{n \geq 1} n \widetilde{E}_\nu \left| \frac{1}{n} \sum_{k=1}^n \widetilde{P}_{\sigma_k}(Y = n) - \widetilde{P}_\nu(Y = n) \right| < \infty.$$

*The $x \log x$ condition is then necessary.*

**Proof.** The condition in the proposition implies that

$$\sum_{n \geq 1} n \left| \frac{1}{n} \sum_{k=1}^n \widetilde{P}_{\sigma_k}(Y = n) - \widetilde{P}_\nu(Y = n) \right| < \infty \qquad \widetilde{P}_\nu\text{-a.s.},$$

which yields

$$\sum_{k=1}^\infty \widetilde{P}_{\sigma_k}(Y > k) = \sum_{k=1}^\infty \sum_{n > k} \widetilde{P}_{\sigma_k}(Y = n)$$

$$= \sum_{n > 1} n \left( \frac{1}{n} \sum_{k=1}^n \widetilde{P}_{\sigma_k}(Y = n) \right)$$

$$\geq \sum_{n > 1} n \widetilde{P}_\nu(Y = n) - \sum_{n \geq 1} n \left| \frac{1}{n} \sum_{k=1}^n \widetilde{P}_{\sigma_k}(Y = n) - \widetilde{P}_\nu(Y = n) \right|$$

$$= \infty$$

since the first term equals $\widetilde{E}_\nu[Y]$, which is infinite if $\widetilde{E}_\nu[\log^+ \bar{\xi}]$ is infinite, this being the case if the $x \log x$ condition does not hold. The second term is finite by assumption. $\square$

Note that if $\{\sigma_n\}$ is Harris recurrent (Alsmeyer (1994), Durrett (2005)) with stationary distribution $\nu$, then the ergodic theorem yields

$$\frac{1}{n} \sum_{k=1}^n \widetilde{P}_{\sigma_k}(Y = j) \to \widetilde{P}_\nu(Y = j) \qquad \widetilde{P}_\nu\text{-a.s.}$$

for all $j$, so our condition means that this convergence is, in some vague sense, "fast enough".

Another generalization of Borel–Cantelli is the Kochen–Stone lemma that states that if $\sum_n P(A_n) = \infty$, then

$$P(A_n \text{ i.o.}) \geq \limsup_n \frac{\{\sum_{k=1}^n P(A_k)\}^2}{\sum_{1 \leq j,k \leq n} P(A_j \cap A_k)}.$$

We can apply this to prove the following result.

**Proposition 6.4.** *Let $A_n = \{\log^+ \bar{\xi}_n > cn\}$. If the (indicators of the) $A_n$ are pairwise negatively correlated, then the $x \log x$ condition is necessary.*

**Proof.** Because $\widetilde{P}_\nu(A_j \cap A_k) \leq \widetilde{P}_\nu(A_j)\widetilde{P}_\nu(A_k)$, we get

$$\widetilde{P}_\nu(A_n \text{ i.o.}) \geq \limsup_n \frac{\{\sum_{k=1}^n \widetilde{P}_\nu(A_k)\}^2}{\sum_{1 \leq j,k \leq n} \widetilde{P}_\nu(A_j \cap A_k)} \geq 1. \qquad \square$$

## Acknowledgements

The author would like to thank Gerold Alsmeyer and Olle Häggström for fruitful discussions related to the results in Section 6, the former of which took place more than ten years ago.